\documentclass{amsart}

\newtheorem{theorem}{Theorem}
\newcommand{\bt}{\begin{theorem}}
\newcommand{\et}{\end{theorem}}
\newtheorem{lemma}{Lemma}
\newcommand{\bl}{\begin{lemma}}
\newcommand{\el}{\end{lemma}}
\newtheorem{corollary}{Corollary}
\newcommand{\bc}{\begin{corollary}}
\newcommand{\ec}{\end{corollary}}
\newcommand{\beq}{\begin{equation}}
\newcommand{\eeq}{\end{equation}}
\newcommand{\benum}{\begin{enumerate}}
\newcommand{\eenum}{\end{enumerate}}
\newcommand{\N}{\ensuremath{ \mathbf N }}
\newcommand{\Z}{\ensuremath{\mathbf Z}}
\newcommand{\Q}{\ensuremath{\mathbf Q}}
\newcommand{\R}{\ensuremath{\mathbf R}}

\newcommand{\mcm}{\ensuremath{ \mathcal M}}

\newcommand{\mct}{\ensuremath{ \mathcal T}}
\DeclareMathOperator{\height}{\text{ht}}

\newcommand{\bmat}{\left(\begin{matrix}}
\newcommand{\emat}{\end{matrix}\right)}
\DeclareMathOperator{\qand}{\quad\text{and}\quad}
\DeclareMathOperator{\qqand}{\qquad\text{and}\qquad}

\usepackage{amsmath,amssymb}
\usepackage[all]{xy}

\title{Free monoids and forests of rational numbers}
\author{Melvyn B. Nathanson}
\address{Department of Mathematics\\
Lehman College (CUNY)\\Bronx, NY 10468}
\email{melvyn.nathanson@lehman.cuny.edu}

\thanks{Supported in part by a grant from the PSC-CUNY Research Award Program.}

\date{\today}

\begin{document}

\begin{abstract}
The Calkin-Wilf tree is an infinite  binary tree whose vertices 
are the positive rational numbers.  
Each such number occurs in the tree exactly once and in the form $a/b$, 
where are $a$ and $b$ are relatively prime positive integers.  
This tree is associated with 
 the matrices $L_1 = \bmat 1 & 0 \\ 1 & 1 \emat$ and 
$R_1 = \bmat 1 & 1 \\ 0 & 1 \emat$, which 
freely generate the monoid $SL_2(\N_0)$ 
of $2 \times 2$ 
matrices with  determinant 1 and nonnegative integral coordinates.  
For other pairs of matrices $L_u$ and $R_v$ that freely generate submonoids 
of $GL_2(\N_0)$, there are forests of infinitely many rooted infinite binary trees 
that partition the set of positive rational numbers, 
and possess a remarkable symmetry property.

\end{abstract}

\maketitle

\section{The Calkin-Wilf tree of rational numbers}

A directed graph is a \emph{rooted infinite  binary tree} 
if it is a tree with the following properties:
\benum
\item[(i)]
Every vertex is the tail of exactly two edges. 
Equivalently, every vertex has \emph{outdegree} 2.
\item[(ii)]
There is a vertex $z$ such that every vertex $v\neq z$ is the head 
of exactly one edge, but $z$ is not the head of any edge.  
Equivalently, every vertex $v\neq z$ has \emph{indegree} 1, and $z$ has indegree 0.
We call $z$ the \emph{root} of the tree.
\item[(iii)]
The graph is connected.
\eenum
In this paper, a \emph{forest} is a directed graph whose connected components are 
rooted infinite  binary trees.

Let $\Q^+$ denote the set of positive rational numbers.  
We call the  rational number $a/b$ \emph{reduced} if $b\geq 1$ and 
the integers $a$ and $b$ are relatively prime.
The Calkin-Wilf tree~\cite{calk-wilf00} is a rooted infinite  binary tree whose vertex set 
is the set of positive reduced rational numbers, 
and whose root is 1.  
In this tree, every positive reduced rational number 
$a/b$ is the tail of two edges.  
The heads of these edges are the positive rational numbers $a/(a+b)$ and $(a+b)/b$.  
We draw this as follows:
\[
\xymatrix{
& \frac{a}{b} \ar[dl]  \ar[dr] & \\
\frac{a}{a+b} &  & \frac{a+b}{b} \\
}
\]
with $a/(a+b)$ on the left and $(a+b)/b$ on the right.  
Note that 
\[
0 < \frac{a}{a+b}  < 1 < \frac{a+b}{b}.
\]
Equivalently, if $w = a/b$, then the generation rule of the tree is 
\beq   \label{Forest:parent-z}
\xymatrix{
&w \ar[dl]  \ar[dr] & \\
\frac{w}{w+1} &  & w+1 \\
}
\eeq
Calkin and Wilf~\cite{calk-wilf00} introduced this enumeration 
of the positive rationals in 2000.  
It is related to the Stern-Brocot sequence~\cite{broc60,ster58}, 
discussed in~\cite{gram-knut-pata94}, 
and has stimulated much recent 
research (e.g. 
\cite{adam10, bate-bund-togn10, bate-mans11,dilc-stol07,glas11,mall11,mans-shat11,nath14f,rezn90}).  
For work related to this paper, see \cite{HMST14,HMST15}.

The first four rows the Calkin-Wilf tree are as follows:
\[
\xymatrix@=0.2cm{
& & & & \frac{1}{1} \ar[dll]  \ar[drr] & & & & \\
& & \frac{1}{2} \ar[dl]  \ar[dr]  &  & & & \frac{2}{1} \ar[dl]  \ar[dr]  & & \\
& \frac{1}{3}   \ar[dl]  \ar[dr]  & & \frac{3}{2}  \ar[dl]  \ar[dr] &  & \frac{2}{3}  \ar[dl]  \ar[dr] &  & \frac{3}{1}  \ar[dl]  \ar[dr] & \\
\frac{1}{4}&  & \frac{4}{3} \quad \frac{3}{5} &  & \frac{5}{2} \quad  \frac{2}{5} &  & 
\frac{5}{3} \quad  \frac{3}{4} & & \frac{4}{1}  \\
}
\]
We enumerate the numbers on the rows of the Calkin-Wilf tree as follows.
Row 0 contains only the number 1.  
Row 1 contains  the numbers 1/2 and 2.
For every nonnegative integer $n$, the $n$th row of the Calkin-Wilf tree 
contains $2^n$ positive reduced rational numbers. 
The $n$th row of the tree is also called the \emph{$n$th generation} of the tree.
We denote the ordered sequence of elements of the $n$th row, 
from left to right, by
$c(n,1), c(n,2), \ldots, c(n,2^n)$.  
For example, $c(2,3) = 2/3$ and $c(3,6) = 5/3$.  
Note that $0 < c(n,2i-1) < 1 < c(n,2i)$ for $i=1,2,\ldots, 2^{n-1}$.

Here are four properties of the Calkin-Wilf tree:
\benum
\item[(i)]
Symmetry formula: 
For every nonnegative integer $n$ and for $i =1,\ldots, 2^n$, 
\[
c(n,i)  c(n,2^n  + 1 - i) = 1.
\]
The proof is by induction on $n$.

\item[(ii)]
Denominator-numerator formula: 
For every positive integer $n$, we have $c(n,1) = 1/(n+1)$ 
and  $c(n,2^n) = n+1$.   
For $j=1,\ldots, 2^n-1$, if $c(n,j) = p/q$, then $c(n,j+1) = q/r$.  
Thus, as we move through the Calkin-Wilf tree from row to row, 
and from left to right across each row, the denominator of each fraction 
in the tree is the numerator of the next fraction in the tree.
This is in Calkin-Wilf~\cite{calk-wilf00}.

\item[(iii)]
Successor formula: 
For every positive integer $n$ and for $j=1,\ldots, 2^n -1$, we have
\[
c(n,j+1) = \frac{1}{    2[c(n,j)]+1-c(n,j)  }
\]
where $[x]$ denotes the integer part of the real number $x$.  
This result is due to Moshe Newman~\cite{aign-zieg04,newm03}.

\item[(iv)]
Row formula:
Let $a/b$ be a positive reduced rational number.
If 
\begin{align*}
\frac{a}{b} 
& = a_0 + \cfrac{1}{a_1 + \cfrac{1}{ a_2 + \cdots + \cfrac{1}{ a_{k-1}
+ \cfrac{1}{ a_k}}}} \\
& = [a_0, a_1,\ldots, a_{k-1}, a_k ]
\end{align*}
is the finite continued fraction of $a/b$, 
then $a/b$ appears on the $n$th row 
of the Calkin-Wilf tree, where $ n = a_0 + a_1 + \cdots + a_{k-1}  + a_k - 1$.
This is discussed in Gibbons, Lester, and Bird~\cite{gibb-lest-bird06}.
\eenum

\section{Freely generated monoids and a symmetry of trees}  

A monoid is a semigroup with an identity.  
Let $GL_2(\R_{\geq 0})$ denote the multiplicative monoid of $2 \times 2$ 
matrices with nonzero determinant and with coordinates in 
the set $\R_{\geq 0}$ of nonnegative real numbers.  
To every matrix 
\[
A =  \bmat a_{1,1} & a_{1,2} \\ a_{2,1} & a_{2,2} \emat \in GL_2(\R_{\geq 0})
\]
we associate the  linear fractional transformation
\[
A(w) = \frac{a_{1,1}w + a_{1,2} }{a_{2,1} w+ a_{2,2}} .
\]
%defined on the set of positive real numbers $t$.  
This is a monoid isomorphism from $GL_2(\R_{\geq 0})$ 
to the monoid of linear fractional transformations with 
nonnegative real coordinates, nonzero determinant, and the binary operation 
of composition of functions.     

The monoid $\mcm(A,B)$ generated by a pair of matrices 
$\{A,B\}$ in $GL_2(\R_{\geq 0})$
consists of all matrices that can be represented as 
products of nonnegative powers of $A$ and $B$.
The matrices $A$ and $B$ \emph{freely generate} this monoid 
if every matrix in $\mcm(A,B)$ has a unique representation 
as a product of powers of $A$ and $B$.

It is well-known (often described as a ``folk theorem'')  that the matrices 
\[
L_1 = \bmat 1 & 0 \\ 1 & 1 \emat
\qqand
R_1 = \bmat 1 & 1 \\ 0 & 1 \emat
\]
freely generate the monoid $SL_2(\N_0)$ 
of $2 \times 2$ 
matrices with  determinant 1 and nonnegative integral coordinates.  
The corresponding linear fractional transformations are  
\[
L_1(w) = \frac{w}{w+1}
\qqand
R_1(w) = w+1.
\]
We observe that  
\[
0 < L_1(w) < 1 < R_1(w)
\]
for all $w \in \Q^+$.  
We can rewrite the generation rule~\eqref{Forest:parent-z}  
of the Calkin-Wilf tree in the form 
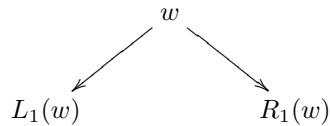
\beq   \label{Forest:parent-L1R1}
\xymatrix{
&w \ar[dl]  \ar[dr] & \\
L_1(w) &  & R_1(w) \\
}
\eeq
That the Calkin-Wilf graph  with vertex set $\Q^+$ is a  tree implies
that  the matrices $L_1$ and $R_1$ freely generate the monoid $\mcm(L_1,R_1)$. 
 
A standard generating set for the group $SL_2(\Z)$ is 
$\left\{  \bmat 1 & 1 \\ 0 & 1 \emat, \bmat 0 & -1 \\ 1 & 0 \emat \right\}$.  
Because 
\[
L_1R_1^{-1}L_1 =  \bmat 0 & -1 \\ 1 & 0 \emat 
\qand
\left( L_1R_1^{-1}L_1\right)^4 =  \bmat 1 & 0 \\ 0 & 1 \emat 
\]
it follows that $\left\{ L_1, R_1\right\}$ generates but does not freely generate $SL_2(\Z)$.

Let $L$ and $R$ be matrices in $GL_2(\R_{\geq 0})$ such that 
\beq          \label{Forest:inequalityLR}
0 < L(w) < 1 < R(w)
\eeq
for all $w \in \Q^+$.  
If the coordinates of $L$ and $R$ are nonnegative integers, then $L(w) \in \Q^+$ 
and  $R(w) \in \Q^+$ for all  $w \in \Q^+$.
For every positive rational number $z$, we can construct 
inductively a directed graph with root $z$ 
such that every vertex is the tail of two edges:   
\beq    \label{Forest:generateLR}
\xymatrix{
&w \ar[dl]  \ar[dr] & \\
L(w) &  & R(w) \\
}
\eeq
Inequality~\eqref{Forest:inequalityLR} and the invertibility of the matrices $L$ and $R$ 
imply that this graph is a rooted infinite binary tree.

A standard application of the ping-pong lemma 
(e.g. Lyndon and Schupp~\cite[pp. 167--168]{lynd-schu77}) 
proves that, for every pair $(u,v)$ of integers with $u \geq 2$ and $v \geq 2$, 
the matrices 
\[
L_u = \bmat 1 & 0 \\ u & 1 \emat
\qqand
R_v = \bmat 1 & v \\ 0 & 1 \emat
\]
generate a free group of rank 2.    
In particular, the nonnegative powers of $L_u$ and $R_v$ generate 
a free monoid.  
The case $u=v=2$ is Sanov's theorem~\cite{sano47a}.  

These are special cases of the following result.

\bt [Nathanson~\cite{nath14g} ]           \label{CWFree:theorem}
Let $A = \bmat a_{1,1} & a_{1,2} \\ a_{2,1} & a_{2,2} \emat $ 
and  $B = \bmat b_{1,1} & b_{1,2} \\ b_{2,1} & b_{2,2} \emat $ 
be matrices in $GL_2(\R_{\geq 0})$.  
If 
\[
a_{1,1} \leq a_{2,1} \qqand  a_{1,2} \leq a_{2,2}
\]
and if 
\[
b_{1,1} \geq b_{2,1}  \qqand  b_{1,2} \geq b_{2,2}
\]
then 
\benum
\item[(i)]
for all $w \in \Q^+$,
\[
0 < A(w) < 1 < B(w)
\]  
\item[(ii)]
the submonoid of $GL_2(\R_{\geq 0})$ 
generated by $A$ and $B$ is free, 
\item[(iii)]
the matrices $A$ and $B$  freely generate $\mcm(A,B)$.
 \eenum
\et

Theorem~\ref{CWFree:theorem} implies that if $u$ and $v$ are are positive 
integers, then the matrices 
\beq      \label{Forest:LuRv}
L_u = \bmat 1 & 0 \\ u & 1 \emat \qand R_v = \bmat 1 & v \\ 0 & 1 \emat
\eeq
freely generate a submonoid of the multiplicative monoid $GL_2(\R_{\geq 0})$, 
and the directed graph $\mct^{(u,v)}_{z}$ with root $z$ and generation rule 
\[
\xymatrix{
& w \ar[dr] \ar[dl]  & \\
L_u(w) = \frac{w}{uw+1} & & R_v(w) = w + v 
}
\]
is a rooted infinite binary tree.  If $a/b$ is a positive reduced fraction and $w = a/b$, 
then the generation rule is 
\beq      \label{Forest:generationRuleabuv}
\xymatrix{
& \frac{a}{b} \ar[dr] \ar[dl]  & \\
\frac{a}{ua+b} & & \frac{a+vb}{b} 
}
\eeq

Let $a/b \in \Q^+$.   
If $a > vb$, then $a/b = R_v((a-vb)/b)$.  
If $b > ua$, then $a/b = L_u(a/(b-ua))$.  
If 
\beq    \label{Forest:orphanIneq}
\frac{1}{u}  \leq \frac{a}{b} \leq v
\eeq
then $a/b$ is an \emph{orphan}, that is, 
$a/b \neq L_u(w)$ and $a/b \neq R_v(w)$ for all $w \in \Q^+$.  
Thus, if $a/b \in \Q^+$ satisfies inequality~\eqref{Forest:orphanIneq}, then 
$a/b$ is a vertex  in a rooted 
infinite binary tree with generation rule~\eqref{Forest:generationRuleabuv} 
if and only if it is the root of the tree.

We define the \emph{height} of the reduced rational number $a/b$ 
by $\height(a/b) = \max\{|a|,|b|\}$.  
If $u$ and $v$ are positive integers and $a/b > 0$, then  
\[
\height(a/b) < a+b \leq ua+b = \height(L_u(a/b)
\]
and 
\[
\height(a/b) < a+b \leq a +vb  = \height(R_v(a/b).
\]
Because the height of every reduced rational number is a positive integer,
it follows that every $a/b \in \Q^+$ has only finitely many ancestors, and so 
every  positive rational number is a vertex in some rooted infinite binary tree 
whose root is a rational number satisfying inequality~\eqref{Forest:orphanIneq}.  
This proves that the forest of such rooted infinite binary trees partitions $\Q^+$.

\emph{Notation:}  For $n = 0,1,2,\ldots$ and $i=1,2,\ldots, 2^n$, we denote  by $c^{(u,v)}_z(n,i)$ 
the $i$th number on the $n$th row of the rooted infinite binary tree with root $z$:
{\small 
\[
\xymatrix@=0.1cm{
& c^{(u,v)}_z(n,i) \ar[dddr] \ar[dddl]  & \\
& & \\
& & \\
 c^{(u,v)}_z(n+1,2i-1) = L_u\left( c^{(u,v)}_z(n,i)  \right) & & 
c^{(u,v)}_z(n+1,2i)  = R_v\left( c^{(u,v)}_z(n,i)  \right)  
}
\]
}

We examine some trees associated with pairs $(u,v)$ of positive integers.  
For $(u,v) = (1,1)$, the unique orphan is $z=1$, 
and we obtain the Calkin-Wilf tree, whose vertex set 
is the set of all positive rational numbers, and $c^{(1,1)}_1(n,i) = c(n,i)$.

Consider the case $(u,v) = (2,2)$.  
In the forest of trees of positive fractions generated by the matrices 
$\bmat 1 & 0 \\ 2 & 1 \emat$ and 
$\bmat 1 & 2 \\ 0 & 1 \emat$, the roots of the trees are the rational numbers $z$ such that 
$1/2 \leq z \leq 2$.  
We consider the trees with roots 1, 3/2, and 2/3.  
For simplicity, we omit the arrows connecting vertices on successive rows.

The first five rows of the tree with root 1 are
{\small
{\[\displaystyle 1\]}
{\[\displaystyle  \frac{1}{3}\quad3\]}
{\[\displaystyle \frac{1}{5}\quad\frac{7}{3}\quad\frac{3}{7}\quad5\]}
{\[\displaystyle \frac{1}{7}\quad{\frac {11}{5}}\quad{\frac {7}{17}}\quad  \frac{13}{3} \quad  \frac{3}{13} \quad{\frac {17}{7}}\quad{\frac {5}{11}}\quad7\]}
{\[\displaystyle \frac{1}{9} \quad{\frac {15}{7}}\quad{\frac {11}{27}}\quad{\frac {21}{5}}\quad{\frac {7}{31}}\quad{\frac {41}{17}}\quad{\frac {13}{29}}\quad{\frac {19}{3}}\quad{\frac {3}{19}}\quad{\frac {29}{13}}\quad{\frac {17}{41}}\quad{\frac {31}{7}}\quad{\frac {5}{21}}\quad{\frac {27}{11}}\quad{\frac {7}{15}} \quad9\]}
}\\
Observe the symmetry in each line:
\beq    \label{Forest:symmetry-uu-1}
c^{(2,2)}_1(n,i) \ c^{(2,2)}_1(n,2^n+1-i) = 1
\eeq
for $i=1,2,\ldots, 2^n$.

The first five rows of the tree with root 3/2 are
{\small 
{\[\displaystyle \frac{3}{2}\]}
{\[\displaystyle \frac{3}{8}\quad \frac{7}{2}\]}
{\[\displaystyle \frac{3}{14}\quad{\frac {19}{8}}\quad{\frac {7}{16}}\quad\frac{11}{2}  \]}
{\[\displaystyle {\frac {3}{20}} \quad {\frac {31}{14}}\quad{\frac {19}{46}}\quad{\frac {35}{8}}\quad{\frac {7}{30}}\quad{\frac {39}{16}}\quad{\frac {11}{24}}\quad\frac{15}{2}\]}
{\[\displaystyle {\frac {3}{26}}\quad{\frac {43}{20}}\quad{\frac {31}{76}}\quad{\frac {59}{14}}\quad{\frac {19}{84}}\quad{\frac {111}{46}}\quad{\frac {35}{78}}\quad{\frac {51}{8}}\quad{\frac {7}{44}}\quad{\frac {67}{30}}\quad{\frac {39}{94}}\quad{\frac {71}{16}}\quad{\frac {11}{46}}\quad{\frac {59}{24}}\quad{\frac {15}{32}}\quad\frac{19}{2}\]} 
}\\
In this case, the symmetry of type~\eqref{Forest:symmetry-uu-1} in each line disappears.  
However, look at 
the first five rows of the tree with the reciprocal root 2/3.
{\small
{\[\displaystyle \frac{2}{3}  \]}
{\[\displaystyle \frac{2}{7}\quad\frac{8}{3}\]}
{\[\displaystyle \frac{2}{11}\quad{\frac {16}{7}}\quad{\frac {8}{19}}\quad\frac{14}{3}\]}
{\[\displaystyle \frac{2}{15}\quad{\frac {24}{11}}\quad{\frac {16}{39}}\quad{\frac {30}{7}}\quad{\frac {8}{35}}\quad{\frac {46}{19}}\quad{\frac {14}{31}}\quad{\frac {20}{3}}\]}
{\[\displaystyle \frac{2}{19}\quad{\frac {32}{15}}\quad{\frac {24}{59}}\quad{\frac {46}{11}}\quad{\frac {16}{71}}\quad{\frac {94}{39}}\quad{\frac {30}{67}}\quad{\frac {44}{7}}\quad{\frac {8}{51}}\quad{\frac {78}{35}}\quad{\frac {46}{111}}\quad{\frac {84}{19}}\quad{\frac {14}{59}}\quad{\frac {76}{31}}\quad{\frac {20}{43}}\quad{\frac {26}{3}}\]} 
}\\
We observe a new symmetry between corresponding lines of the two trees:   
\beq    \label{Forest:symmetry-uu}
c^{(2,2)}_{3/2}(n,i) \ c^{(2,2)}_{2/3}(n,2^n+1-i) = 1
\eeq
for $i =1,2,\ldots, 2^n$.  

Consider next the case $(u,v) = (5,4)$ and the forest of positive fractions 
generated by the matrices 
$\bmat 1 & 0 \\ 5 & 1 \emat$ and 
$\bmat 1 & 4 \\ 0 & 1 \emat$.   
The first five rows of the tree with root 3/2 are
{\small 
{\[\displaystyle {\frac {3}{2}}  \]}
{\[\displaystyle {\frac {3}{17}} \quad {\frac {11}{2}}   \]}
{\[\displaystyle {\frac {3}{32}}\quad{\frac {71}{17}}\quad{\frac {11}{57}}\quad {\frac {19}{2}}   \]}
{\[\displaystyle {\frac {3}{47}}\quad{\frac {131}{32}}\quad{\frac {71}{372}}\quad{\frac {139}{17}}\quad{\frac {11}{112}}\quad{\frac {239}{57}}\quad{\frac {19}{97}}\quad{\frac {27}{2}}\]}
{\[\displaystyle {\frac {3}{62}}\quad{\frac {191}{47}}\quad{\frac {131}{687}}\quad{\frac {259}{32}}\quad{\frac {71}{727}}\quad{\frac {1559}{372}}\quad{\frac {139}{712}}\quad{\frac {207}{17}}\quad{\frac {11}{167}}\quad{\frac {459}{112}}\quad{\frac {239}{1252}}\quad{\frac {467}{57}}\quad{\frac {19}{192}}\quad{\frac {407}{97}}\quad{\frac {27}{137}}\quad{\frac {35}{2}}\]} 
}\\
Again there is no symmetry in each line.

We look at the first five rows of the tree with the reciprocal root 2/3.
{\small 
{\[\displaystyle    {\frac {2}{3}} \]}
{\[\displaystyle   {\frac {2}{13}}\quad  {\frac {14}{3}}\]}
{\[\displaystyle   {\frac {2}{23}}\quad{\frac {54}{13}}\quad{\frac {14}{73}}\quad{\frac {26}{3}}\]}
{\[\displaystyle {\frac {2}{33}}\quad{\frac {94}{23}}\quad{\frac {54}{283}}\quad{\frac {106}{13}}\quad{\frac {14}{143}}\quad{\frac {306}{73}}\quad{\frac {26}{133}}\quad{\frac {38}{3}}\]}
{\[\displaystyle {\frac {2}{43}}\quad{\frac {134}{33}}\quad{\frac {94}{493}}\quad{\frac {186}{23}}\quad{\frac {54}{553}}\quad{\frac {1186}{283}}\quad{\frac {106}{543}}\quad{\frac {158}{13}}\quad{\frac {14}{213}}\quad{\frac {586}{143}}\quad{\frac {306}{1603}}\quad{\frac {598}{73}}\quad{\frac {26}{263}}\quad{\frac {558}{133}}\quad{\frac {38}{193}}\quad{\frac {50}{3}}\]}
}

In contrast to the case $(u,v) = (2,2)$, when $(u,v) = (5,4)$, 
we do not observe a symmetry of the form~\eqref{Forest:symmetry-uu} 
between corresponding lines of the trees with reciprocal roots 3/2 and 2/3, 
However, consider the first five rows of the tree with $(u,v) = (4,5)$,
generated by the matrices 
$\bmat 1 & 0 \\ 4 & 1 \emat$ and 
$\bmat 1 & 5 \\ 0 & 1 \emat$, and with the root 2/3: 
\\
{\small
{\[\displaystyle  {\frac {2}{3}} \]}
{\[\displaystyle   {\frac {2}{11}}\quad{\frac {17}{3}}\]}
{\[\displaystyle   {\frac {2}{19}}\quad{\frac {57}{11}}\quad{\frac {17}{71}}\quad{\frac {32}{3}}\]}
{\[\displaystyle {\frac {2}{27}}\quad{\frac {97}{19}}\quad{\frac {57}{239}}\quad{\frac {112}{11}}\quad{\frac {17}{139}}\quad{\frac {372}{71}}\quad{\frac {32}{131}}\quad{\frac {47}{3}}\]}
{\[\displaystyle {\frac {2}{35}}\quad{\frac {137}{27}}\quad{\frac {97}{407}}\quad{\frac {192}{19}}\quad{\frac {57}{467}}\quad{\frac {1252}{239}}\quad{\frac {112}{459}}\quad{\frac {167}{11}}\quad{\frac {17}{207}}\quad{\frac {712}{139}}\quad{\frac {372}{1559}}\quad{\frac {727}{71}}\quad{\frac {32}{259}}\quad{\frac {687}{131}}\quad{\frac {47}{191}}\quad{\frac {62}{3}}\]} 
}\\
A beautiful symmetry reappears:  
In the tree constructed from the pair $(5,4)$ with the root 3/2, 
and in the tree constructed from the reversed pair $(4,5)$ with the root 2/3, 
we find 
\beq    \label{Forest:symmetry-uv}
c^{(5,4)}_{3/2}(n,i) \ c^{(4,5)}_{2/3}(n,2^n+1-i) = 1
\eeq
for $i=1,2,\ldots, 2^n$.  
We shall prove that this identity holds for all pairs $(u,v)$ of real numbers 
such that $u \geq 1$ and $v \geq 1$, and for all roots $z$.

\section{Proof of symmetry}

\bt [Symmetry]        \label{CWPingPong:theorem:c-symmetry}
Let $z$ be a variable, and let $u$ and $v$ be positive integers.  
For all $n \in \N_0$ and $i =1,2,\ldots, 2^n$,  
\[
c^{(u,v)}_z(n,i) c^{(v,u)}_{z^{-1}}(n,2^n+1-i) = 1.
\] 
If $u = v \geq 1$, then for all $n \in \N_0$ and $i=1,2,\ldots, 2^n$,  
\[
c^{(u,u)}(n,i) c^{(u,u)}(n,2^n+1-i) = 1.
\] 
\et

If $u=v=1$, then this is the familiar symmetry of the Calkin-Wilf tree.

\begin{proof}
The proof is by induction on the row number $n$.
For $n=0$ and $i=1$, we have $c^{(u,v)}_z(0,1) = z$ and so 
\[
c^{(u,v)}_z(n,i) c^{(v,u)}_{z^{-1}}(n,2^n+1-i) = z z^{-1} = 1.
\]
Let $n \geq 0$, and assume that the Theorem holds for row $n$.  
For $i = 1,2,\ldots, 2^n$, we have 
\[
c^{(u,v)}_{z}(n,i) = \frac{1}{  c^{(v,u)}_{z^{-1}}(n,2^n+1-i) }.
\]
It follows that 
\begin{align*}
c^{(u,v)}_z(n+1,2i-1) 
& =  L_u\left( c^{(u,v)}_{z}(n,i)  \right) \\
& = \frac{  c^{(u,v)}_{z}(n,i) }{  uc^{(u,v)}_{z}(n,i) +1  } \\
& = \frac{ 1  }{  c^{(v,u)}_{z^{-1}}(n,2^n+1-i) +u } \\
& = \frac{ 1  }{ R_u\left(  c^{(v,u)}_{z^{-1}}(n, 2^n + 1 - i  )  \right) } \\
& = \frac{ 1  }{ c^{(v,u)}_{z^{-1}}(n+1, 2\left(  2^n + 1 - i \right) ) } \\
& = \frac{ 1  }{ c^{(v,u)}_{z^{-1}}(n+1,2^{n+1} + 1 - (2i-1))}.
\end{align*}
Similarly,
\begin{align*}
c^{(u,v)}_z(n+1,2i) 
& =  R_v\left( c^{(u,v)}_{z}(n,i)  \right) \\
& =  c^{(u,v)}_{z}(n,i)  + v \\
& = \frac{ v c^{(v,u)}_{z^{-1}}(n,2^n+1-i) +1 }{  c^{(v,u)}_{z^{-1}}(n,2^n+1-i)  } \\
& = \frac{1}{  L_v( c^{(v,u)}_{z^{-1}}(n,2^n+1-i) )  } \\
& = \frac{1}{ c^{(v,u)}_{z^{-1}}(n+1, 2\left(2^n+1-i \right) -1)  } \\
& =  \frac{1}{ c^{(v,u)}_{z^{-1}}(n+1, 2^{n+1}+1- 2i)  }.
\end{align*}
This completes the proof.
\end{proof}

\section{Open problems}

\benum

\item
If $u$ and $v$ are are positive integers and $L_u$ and $R_v$ are the matrices 
defined by equation~\eqref{Forest:LuRv}, then every positive rational number  $w$ is a vertex 
in a tree in the forest of rooted infinite binary trees generated by $L_u$ and $R_v$.  
In particular, the positive rational number $w$ has only finitely many ancestors in this tree.   
It is not known if there are matrices $L$ and $R$  in $GL_2(\R_{\geq 0})$ 
(that may or may not freely generate the monoid $\mcm(L,R)$) 
such that there exists a positive rational number $z$ that has \emph{infinitely many} ancestors
in the directed graph with vertex $z$ and generation rule~\eqref{Forest:generateLR}.

\item
Do there exist pairs of matrices $L$ and $R$ that do not 
satisfy inequality~\eqref{Forest:inequalityLR} but establish a partition of the  
positive rational numbers into pairwise disjoint rooted infinite binary trees?

\item
Find analogues of properties~(ii),~(iii), and~(iv) of the Calkin-Wilf tree 
that apply to the trees $\mct^{(u,v)}_z$.  (For the trees of linear fractional 
transformations associated to the pair $(1,1)$, see Nathanson~\cite{nath14f}.

\item
Let $m \geq 3$.   Do there exist $m$ matrices $A_1, \ldots, A_m$ in $GL_2(\N_))$ such that 
the $m$-ary generation rule 
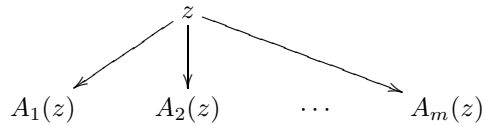
\[
\xymatrix{
 & z \ar[drr] \ar[d]   \ar[dl]  & &\\  
A_1(z) & A_2(z) & \cdots & A_m(z)
}
\]
determines a forest of rooted infinite $m$-ary trees that partition the positive rational numbers?

\eenum

\def\cprime{$'$} \def\cprime{$'$} \def\cprime{$'$} \def\cprime{$'$}
\providecommand{\bysame}{\leavevmode\hbox to3em{\hrulefill}\thinspace}
\providecommand{\MR}{\relax\ifhmode\unskip\space\fi MR }
\providecommand{\MRhref}[2]{%
  \href{http://www.ams.org/mathscinet-getitem?mr=#1}{#2}
}
\providecommand{\href}[2]{#2}

\end{document}